\documentclass[a4paper,12pt]{article}
\usepackage{amssymb,amsmath,amsthm}
\usepackage{indentfirst}
\usepackage{amssymb,amsmath,amsthm}
\usepackage{euler}
\usepackage{hyperref,hypbmsec}
\newtheorem{Th}{\hskip\parindent Theorem}
\newtheorem{Le}{\hskip\parindent Lemma}

\newtheorem{Sl}{\hskip\parindent Corollary}

\newcounter{propet}

\renewcommand{\le}{\leqslant}\renewcommand{\ge}{\geqslant}

\begin{document}
\author{D.\,A.\,Frolenkov\footnote{The research was supported by the grant RFBR № 11-01-00759-a} }
\title{
On a problem of Dobrowolski--Williams.}
\date{}
\maketitle
\begin{abstract}
In this paper we prove  new upper bounds for the sum $\sum_{n=a+1}^{a+N}f(n)$, for a certain class of arithmetic functions $f$. Our results improve the previous results of G. Bachman and $\mbox{L. Rachakonda.}$
\end{abstract}
\setcounter{Zam}0
\section{Introduction}
For any positive real numbers $A$, $B$ and a natural number $q$, denote by $\textit{F}=\textit{F}_{A,B}(q)$ the class of all functions $f:\mathbb{Z}\rightarrow\mathbb{C}$ satisfying the conditions
\begin{gather}
|f(n)|\le A \quad\mbox{for all}\quad n\in\mathbb{Z},
\end{gather}
\begin{gather}
f(n+q)=f(n) \quad\mbox{for all} \quad n\in\mathbb{Z},
\end{gather}
\begin{gather}
\sum\limits_{n=1}^q\left|\sum\limits_{k=1}^Kf(n+k)\right|^2\le BqK  \quad\mbox{for all natural numbers} \quad K.
\end{gather}
Dobrowolski and Williams \cite{Dobro-Williams} proved that the estimate
\begin{gather}\label{DW}
\left|\sum_{n=a+1}^{a+N}f(n)\right|\le\frac{\sqrt{B}}{2\log2}\sqrt{q}\log q+3A\sqrt{q}
\end{gather}
holds for all  $f\in\textit{F}_{A,B}(q).$ Bachman and Rachakonda \cite{Bachman-Racha} improved their result and obtained that
\begin{gather}\label{BR}
\left|\sum_{n=a+1}^{a+N}f(n)\right|\le\frac{\sqrt{B}}{3\log3}\sqrt{q}\log q+\left(5\sqrt{B}+\frac{3}{2}A\right)\sqrt{q}.
\end{gather}
In this paper we improve \eqref{BR}.\par
Let $\{q_n\}_{n=-2}^{\infty}$ be a sequence of integers such that
\begin{gather}\label{q}
q_{-2}=1,\,q_{-1}=1,\,q_n=2q_{n-1}+q_{n-2}\quad \mbox{for}\, n\ge 0.
\end{gather}
Then
\begin{gather}\label{qn}
q_n=\left(\frac{3}{2}+\sqrt{2}\right)\lambda_1^n+\left(\frac{3}{2}-\sqrt{2}\right)\lambda_2^n
\end{gather}
where
\begin{gather}
\lambda_1=1+\sqrt{2},\,\lambda_2=1-\sqrt{2}.
\end{gather}
Let $\{p_n\}_{n=-2}^{\infty}$ be a sequence of integers such that
\begin{gather}\label{p}
p_{-2}=0,\, p_{-1}=0,\, p_n=2p_{n-1}+p_{n-2}+\frac{q_{n-1}+q_{n-2}}{2} \quad\mbox{for}\, n\ge 0.
\end{gather}
Then
\begin{gather}\label{pn}
p_n=\left(\frac{1}{2}+\frac{5\sqrt{2}}{16}\right)\lambda_1^n+\left(\frac{1}{2}-\frac{5\sqrt{2}}{16}\right)\lambda_2^n+
\frac{n}{8}\lambda_1^{n+2}+\frac{n}{8}\lambda_2^{n+2}.
\end{gather}
Define the quantity $\delta_n$,$n\ge0$ by
\begin{gather}\label{deltan}
\delta_n=\frac{p_n}{q_n\log q_n}.
\end{gather}
\begin{Th}\label{Th1}
For any $n\ge0$, $q\ge q_n^6$ and any $f\in\textit{F}_{A,B}(q)$, we have
\begin{gather*}
\left|\sum_{m=a+1}^{a+N}f(m)\right|\le\sqrt{B}\delta_n\sqrt{q}\log q+
\left(\sqrt{B}\sqrt{q_n+\frac{1}{q_n^{2}}}+\sqrt{B}(\sqrt{2}-1)\frac{2p_n}{q_n-1}+\frac{1}{2}A\right)\sqrt{q}+\psi_n(q,B),
\end{gather*}
where
\begin{gather*}
\psi_n(q,B)=\sqrt{B}\left(q_n+\frac{1}{q_n^{2}}\right)\left(\delta_n\log q+(\sqrt{2}-1)\frac{2p_n}{q_n-1}\right).
\end{gather*}
\end{Th}
The aim of the Theorem \ref{Th1}  is to improve the constant $\delta_0=\frac{1}{3\log 3}$ in \eqref{BR} (see Table \ref{tab1}).
\begin{center}
\begin{tabular}{|c|c|c|c|}
\hline
$n$  & $q_n$  & $p_n$ & $\delta_n$\\\hline
0 & 3  & 1  & 0.303413 \\\hline
1 & 7  & 4  & 0.293656 \\\hline
2 & 17 & 14 & 0.290670 \\\hline
3 & 41 & 44 & 0.288986 \\\hline
4 & 99 & 131 & 0.287965 \\\hline
\multicolumn{4}{c}{\phantom{\Large G} Table 1}\label{tab1}
\end{tabular}
\end{center}
By  \eqref{qn}, \eqref{pn}, \eqref{deltan} one has
\begin{gather*}
\lim_{n\rightarrow\infty}\delta_n=\frac{1}{4\log(1+\sqrt{2})}=0.283676\ldots.
\end{gather*}
To prove Theorem \ref{Th1} we extend the method of Bachman and Rachakonda.
\section{Proof of the theorem}
If $f\in\textit{F}_{A,B}(q)$ then $\sum\limits_{n=1}^qf(n)=0$. So, we may assume that
\begin{gather}\label{Nq2}
N<\frac{q}{2}.
\end{gather}
If $N\le\frac{\sqrt{q}}{2}$ then
$$\left|\sum_{n=a+1}^{a+N}f(n)\right|\le A\frac{\sqrt{q}}{2}$$
and Theorem \ref{Th1} is proved. So, we may assume that
\begin{gather}\label{N}
N\ge\frac{\sqrt{q}}{2}\ge\frac{q_n^3}{2}.
\end{gather}
As in \cite{Bachman-Racha} we define triangular sums $T_{-}(x,y),T_{+}(x,y)$ by
\begin{gather}\label{T-}
T_{-}(x,y)=\sum_{i=x+1}^{x+y}\sum_{j=x+1-i}^0f(i+j)=\sum_{k=1}^y(y+1-k)f(x+k).
\end{gather}
\begin{gather}\label{T+}
T_{+}(x,y)=\sum_{i=x+1}^{x+y}\sum_{j=0}^{x+y-i}f(i+j)=\sum_{k=1}^ykf(x+k)
\end{gather}
and the square sum $S(x,y)$ by
\begin{gather}\label{S}
S(x,y)=\sum_{i=x+1}^{x+y}\sum_{j=0}^{y-1}f(i+j).
\end{gather}
Observe that for any integer numbers  $x,y,u,v,k$ one has
\begin{gather}\label{transformation}
\sum_{i=x}^{x+y}\sum_{j=u}^{u+v}f(i+j)=\sum_{i=x-k}^{x+y-k}\,\sum_{j=u+k}^{u+v+k}f(i+j).
\end{gather}
Let $K, 1\le K\le N$. It was shown in \cite[(2.2),(2.3)]{Bachman-Racha} that
\begin{gather}\label{1estimate}
\left|\sum_{n=a+1}^{a+N}f(n)\right|\le\sqrt{B q}\sqrt{\frac{N}{K}}+
\frac{1}{K}\left|-T_{-}(a,K)+T_{+}(a+N,K)\right|.
\end{gather}
Let $\tau$ be the integer satisfying
\begin{gather}\label{NN}
\frac{q_n^{\tau}-1}{2}\le N<\frac{q_n^{\tau+1}-1}{2}.
\end{gather}
By \eqref{N} we have $\tau\ge3.$
Let $\{K_i\}_{i=0}^{\tau}$ be a sequence of integers such that
\begin{gather}\label{Ki}
K_i=\frac{q_n^{\tau-i}-1}{2}
\end{gather}
and $$K=K_0=\frac{q_n^{\tau}-1}{2}.$$
Let $\{c_n\}_{n=-2}^{\infty}$ be a sequence of integers such that
\begin{gather}\label{cn}
c_n=\frac{q_n-1}{2}.
\end{gather}
Note that
\begin{gather}\label{cq}
c_j-c_{j-1}=\frac{q_{j-1}+q_{j-2}}{2}.
\end{gather}
By \eqref{q},  we have
\begin{gather}\label{c}
c_n=2c_{n-1}+c_{n-2}+1.
\end{gather}
Then for $0\le i\le\tau-1$ one has
\begin{gather}\label{K0N}
K_i=q_nK_{i+1}+c_n
\end{gather}
and
\begin{gather}
N<q_nK_{0}+c_n.
\end{gather}
So
\begin{gather}\label{NK0}
\frac{N}{K_0}<q_n+\frac{q_n-1}{q_n^{\tau}-1}<q_n+\frac{1}{q_n^{\tau-1}}<q_n+\frac{1}{q_n^{2}}.
\end{gather}
\begin{Le}\label{recurs-}
For any $1\le i\le \tau$ one has
\begin{gather*}
T_{-}(a,K_{i-1})=2T_{-}(a,q_{n-1}K_i+c_{n-1})-T_{+}(a-q_{n-2}K_i-c_{n-2},q_{n-2}K_i+c_{n-2})+\\+
S(a-q_{n-2}K_i-c_{n-2},(q_{n-2}+q_{n-1})K_i+c_n-c_{n-1}).
\end{gather*}
\begin{proof}
\par
\begin{picture}(250,250)\label{fig}
\put(20,230){
\begin{picture}(250,150)
\put(0,-190){\vector(0,1){220}}
\put(0,0){\vector(1,0){250}}
\put(-10,35){i}
\put(240,10){j}
\put(220,0){\line(0,-1){170}}
\put(50,0){\line(1,-1){170}}
\put(110,0){\line(0,-1){110}}
\put(110,-110){\line(1,0){110}}
\put(105,0){\line(0,-1){55}}
\put(220,-115){\line(-1,0){55}}
\put(110,-65){\line(1,-1){45}}
\put(48,8){A}
\put(222,8){B}
\put(225,-168){C}
\put(115,-12){D}
\put(105,-120){E}
\put(100,8){$D'$}
\put(222,-102){F}
\put(222,-118){$F'$}
\put(92,-60){H}
\put(100,-80){G}
\put(140,-120){M}
\put(160,-132){N}
\put(155,-200){fig. 1}
\end{picture}}
\end{picture}
\par
Put
\begin{gather*}
A(a+1,0);\; B(a+q_{n}K_i+c_{n},0);\; C(a+q_{n}K_i+c_{n},-q_{n}K_i-c_{n}+1).
\end{gather*}
Then
\begin{gather}\label{1}
\sum_{ABC}f(i+j)=T_{-}(a,q_{n}K_i+c_{n})=T_{-}(a,K_{i-1}).
\end{gather}
To prove Lemma \ref{recurs-} we make a partition of the $\triangle ABC$ onto three triangles and the square \mbox{(see fig. \ref{fig}).} So
\begin{gather}\label{0}
\sum_{ABC}f(i+j)=\sum_{AD'H}f(i+j)+\sum_{NF'C}f(i+j)+\sum_{DEFB}f(i+j)-\sum_{EGM}f(i+j).
\end{gather}
Let
\begin{gather}\label{s}
s=(q_n-q_{n-1})K_i+(c_n-c_{n-1}-1).
\end{gather}
By \eqref{q}, \eqref{c} we have
\begin{gather}\label{s1}
s=(q_{n-1}+q_{n-2})K_i+(c_n-c_{n-1}-1)=(q_{n-1}+q_{n-2})K_i+(c_{n-1}+c_{n-2}).
\end{gather}
Put
\begin{gather*}
D(a+q_{n-1}K_i+c_{n-1}+1,0); E(a+q_{n-1}K_i+c_{n-1}+1,-s);
F(a+q_{n}K_i+c_{n},-s).
\end{gather*}
Applying \eqref{transformation}, we have
\begin{gather}\label{2.0}
\sum_{DEFB}f(i+j)=\sum_{D_1E_1F_1B_1}f(i+j),
\end{gather}
where
\begin{gather*}
D_1(a+q_{n-1}K_i+c_{n-1}+1-s,s);\; E_1(a+q_{n-1}K_i+c_{n-1}+1-s,0);\\
F_1(a+q_{n}K_i+c_{n}-s,0);\;B_1(a+q_{n}K_i+c_{n}-s,s).
\end{gather*}
By \eqref{s}, \eqref{s1}
\begin{gather*}
D_1(a-q_{n-2}K_i-c_{n-2}+1,s);\; E_1(a-q_{n-2}K_i-c_{n-2}+1,0);\\
F_1(a+q_{n-1}K_i+c_{n-1}+1,0);\;B_1(a+q_{n-1}K_i+c_{n-1}+1,s).
\end{gather*}
Applying \eqref{S}, \eqref{2.0}, we have
\begin{gather}\label{2}
\sum_{DEFB}f(i+j)=S(a-q_{n-2}K_i-c_{n-2},(q_{n-2}+q_{n-1})K_i+c_n-c_{n-1}).
\end{gather}
Put
\begin{gather*}
  D'(a+q_{n-1}K_i+c_{n-1},0);\;H(a+q_{n-1}K_i+c_{n-1},-q_{n-1}K_i-c_{n-1}+1).
\end{gather*}
By \eqref{T-}, \eqref{transformation}, we have
\begin{gather}\label{3}
\sum_{AD'H}f(i+j)=T_{-}(a,q_{n-1}K_i+c_{n-1}).
\end{gather}
Put
\begin{gather*}
F'(a+q_{n}K_i+c_{n},-s-1);  N(a+s+2,-s-1).
\end{gather*}
Applying \eqref{transformation}, we have
\begin{gather}\label{4.0}
\sum_{NF'C}f(i+j)=\sum_{N_1F_1'C_1}f(i+j),
\end{gather}
where
\begin{gather*}
F_1'(a+q_{n}K_i+c_{n}-s-1,0);\;N_1(a+1,0);\;C_1(a+q_{n}K_i+c_{n}-s-1,-q_{n}K_i-c_{n}+s+2).
\end{gather*}
Note that
$\triangle N_1F_1'C_1=\triangle AD'H$.
By \eqref{4.0},\;\eqref{3}, we have
\begin{gather}\label{4}
\sum_{NF'C}f(i+j)=T_{-}(a,q_{n-1}K_i+c_{n-1}).
\end{gather}
Put
\begin{gather*}
  G(a+q_{n-1}K_i+c_{n-1}+1,-q_{n-1}K_i-c_{n-1}-1);\; M(a+s,-s).
\end{gather*}
Applying \eqref{transformation}, we have
\begin{gather}\label{5.0}
\sum_{EGM}f(i+j)=\sum_{E_1G_1M_1}f(i+j),
\end{gather}
where
\begin{gather*}
E_1(a+q_{n-1}K_i+c_{n-1}+1-s,0);\;G_1(a+q_{n-1}K_i+c_{n-1}+1-s,-q_{n-1}K_i-c_{n-1}-1+s);\;M_1(a,0).
\end{gather*}
By \eqref{T+},\;\eqref{s1}, we have
\begin{gather}\label{5}
\sum_{EGM}f(i+j)=T_{+}(a-q_{n-2}K_i-c_{n-2},q_{n-2}K_i+c_{n-2}).
\end{gather}
Applying \eqref{1}, \eqref{2}, \eqref{3}, \eqref{4}, \eqref{5} to  \eqref{0} one has
\begin{gather*}
T_{-}(a,K_{i-1})=2T_{-}(a,q_{n-1}K_i+c_{n-1})-T_{+}(a-q_{n-2}K_i-c_{n-2},q_{n-2}K_i+c_{n-2})+\\+
S(a-q_{n-2}K_i-c_{n-2},(q_{n-2}+q_{n-1})K_i+c_n-c_{n-1}).
\end{gather*}
Lemma is proved.
\end{proof}
\end{Le}
\begin{Le}\label{recurs+}
For any $1\le i\le \tau$ one has
\begin{gather*}
T_{+}(a-K_{i-1},K_{i-1})=2T_{+}(a-q_{n-1}K_i-c_{n-1},q_{n-1}K_i+c_{n-1})-T_{-}(a,q_{n-2}K_i+c_{n-2})+\\+
S(a-q_{n}K_i-c_{n},(q_{n-2}+q_{n-1})K_i+c_n-c_{n-1}).
\end{gather*}
\begin{proof}
This Lemma can be proved in the same way as Lemma \ref{recurs-}.
\end{proof}
\end{Le}
Let
\begin{gather}
S_{+}^{(n)}(a,K_i)=S(a-q_{n}K_i-c_{n},(q_{n-2}+q_{n-1})K_i+c_n-c_{n-1}),\notag\\
S_{-}^{(n)}(a,K_i)=S(a-q_{n-2}K_i-c_{n-2},(q_{n-2}+q_{n-1})K_i+c_n-c_{n-1}).\label{Sn}
\end{gather}
Let
\begin{gather*}
T_{+}^{(n)}(a,K_i)=T_{+}(a-q_{n}K_i-c_{n},q_{n}K_i+c_{n}),\quad
T_{-}^{(n)}(a,K_i)=T_{-}(a,q_{n}K_i+c_{n}).
\end{gather*}
Note that
\begin{gather}\label{Tn}
T_{-}^{(n)}(a,K_i)=T_{-}(a,K_{i-1}),\qquad  T_{+}^{(n)}(a,K_i)=T_{+}(a-K_{i-1},K_{i-1}).
\end{gather}
Consecutive application of Lemma \ref{recurs-} and  Lemma \ref{recurs+} gives the following result.
\begin{Le}\label{recurs}
For $n\ge0$, $1\le i\le\tau-1$ we have
\begin{gather*}
T_{-}^{(n)}(a,K_{i})=\alpha_nT_{-}(a,K_{i})-\beta_nT_{+}(a,K_{i})+
\sum_{j=0}^na_{j,n}S_{-}^{(j)}(a,K_{i})-\sum_{j=0}^{n-2}b_{j,n}S_{+}^{(j)}(a,K_{i}),
\end{gather*}
\begin{gather*}
T_{+}^{(n)}(a,K_{i})=\alpha_nT_{+}(a,K_{i})-\beta_nT_{-}(a,K_{i})+
\sum_{j=0}^na_{j,n}S_{+}^{(j)}(a,K_{i})-\sum_{j=0}^{n-2}b_{j,n}S_{-}^{(j)}(a,K_{i}),
\end{gather*}
where
\begin{gather*}
\alpha_n=\frac{q_n+1}{2},\quad\beta_n=\frac{q_n-1}{2}
\end{gather*}
and $\{a_{j,n}\}_{j=0}^{n}$, $\{b_{j,n}\}_{j=0}^{n-2}$ are integers.
\begin{proof}
We prove this statement by induction. For $n=0$ the result follows from
\begin{gather}\label{T0-}
T_{-}(a,q_{0}K_i+c_{0})=T_{-}(a,3K_i+1)=2T_{-}(a,K_i)-T_{+}(a-K_i,K_i)+S(a-K_i,2K_i+1),
\end{gather}
\begin{gather}\label{T0+}
T_{+}(a-3K_{i}-1,3K_{i}+1)=2T_{+}(a-K_i,K_i)-T_{-}(a,K_i)+S(a-3K_i-1,2K_i+1),
\end{gather}
(see \cite{Bachman-Racha}). For $n=1$ by Lemma \ref{recurs-}, we have
\begin{gather*}
T_{-}^{(1)}(a,K_{i})=T_{-}(a,K_{i-1})=2T_{-}(a,q_{0}K_i+c_{0})-T_{+}(a-q_{-1}K_i-c_{-1},q_{-1}K_i+c_{-1})+\\+
S(a-q_{-1}K_i-c_{-1},(q_{-1}+q_{0})K_i+c_1-c_{0}).
\end{gather*}
Applying \eqref{T0-}, we obtain
\begin{gather}\label{T1-}
T_{-}^{(1)}(a,K_{i})=4T_{-}(a,K_i)-3T_{+}(a-K_i,K_i)+
S(a-K_i,4K_i+2)+2S(a-K_i,2K_i+1).
\end{gather}
In the same way we can prove that
\begin{gather}\label{T1+}
T_{+}^{(1)}(a,K_{i})=4T_{+}(a-K_i,K_i)-3T_{-}(a,K_i)+
S(a-7K_i-3,4K_i+2)+\notag\\+2S(a-3K_i-1,2K_i+1).\label{T1+}
\end{gather}
If our  formulas are proved for $k\le n-1$ then by Lemma \ref{recurs-} and Lemma \ref{recurs+}, we have
\begin{gather*}
T_{-}^{(n)}(a,K_{i})=2T_{-}^{(n-1)}(a,K_{i})-T_{+}^{(n-2)}(a,K_{i})+S_{-}^{(n)}(a,K_{i})=\\=
(2\alpha_{n-1}+\beta_{n-2})T_{-}^(a,K_{i})-(\alpha_{n-2}+2\beta_{n-1})T_{+}(a,K_{i})+\\+
\Biggl(
S_{-}^{(n)}(a,K_{i})+2a_{n-1,n-1}S_{-}^{(n-1)}(a,K_{i})+2a_{n-2,n-1}S_{-}^{(n-2)}(a,K_{i})+
2a_{n-3,n-1}S_{-}^{(n-3)}(a,K_{i})+\\+\sum_{j=0}^{n-4}(2a_{j,n-1}+b_{j,n-2})S_{-}^{(j)}(a,K_{i})
\Biggl)-\\-
\left(
a_{n-2,n-2}S_{+}^{(n-2)}(a,K_{i})+\sum_{j=0}^{n-3}(a_{j,n-2}+2b_{j,n-1})S_{+}^{(j)}(a,K_{i})
\right).
\end{gather*}
So
\begin{gather*}
\left\{
  \begin{array}{ll}
   \alpha_n=2\alpha_{n-1}+\beta_{n-2} , &  \\
   \beta_n=\alpha_{n-2}+2\beta_{n-1}.
    \end{array}
\right.
\end{gather*}
By \eqref{T0-},\;\eqref{T1-}, we have
\begin{gather*}
\alpha_0=2,\;\beta_0=1,\;\alpha_1=4,\;\beta_1=3.
\end{gather*}
By the definition of $q_n$\eqref{q}, we have
\begin{gather}\label{alphabeta}
\alpha_n=\frac{q_n+1}{2},\quad\beta_n=\frac{q_n-1}{2}.
\end{gather}
For sequences $\{a_{j,n}\}_{j=0}^{n}$, $\{b_{j,n}\}_{j=0}^{n-2}$ we have
\begin{gather}
\left\{
  \begin{array}{lll}
   a_{n,n}=1 , &  \\
   a_{j,n}=2a_{j,n-1}&\quad\mbox{for} \quad n-3\le j\le n-1, &  \\
   a_{j,n}=2a_{j,n-1}+b_{j,n-2}&\quad \mbox{for} \quad0\le j\le n-4, &  \\
   b_{n-2,n}=a_{n-2,n-2}, &  \\
   b_{j,n}=a_{j,n-2}+2b_{j,n-1}& \quad\mbox{for}\quad 0\le j\le n-3.
    \end{array}
\right.\label{sistem}
\end{gather}
\end{proof}
\end{Le}
\begin{Le}\label{Pn}
For $n\ge0$  we have
\begin{gather*}
2p_n=\sum_{j=0}^na_{j,n}(q_{j-1}+q_{j-2})+\sum_{j=0}^{n-2}b_{j,n}(q_{j-1}+q_{j-2}).
\end{gather*}
\begin{proof}
We prove this statement by induction. For $n=0$ we have $2p_0=2$ and by \eqref{sistem},\,\eqref{q} one has
$$a_{0,0}(q_{-1}+q_{-2})=2a_{0,0}=2.$$
It follows from  \eqref{sistem} that
\begin{gather}
\sum_{j=0}^na_{j,n}(q_{j-1}+q_{j-2})+\sum_{j=0}^{n-2}b_{j,n}(q_{j-1}+q_{j-2})=
\sum_{j=0}^{n-1}2a_{j,n-1}(q_{j-1}+q_{j-2})+\sum_{j=0}^{n-3}2b_{j,n-1}(q_{j-1}+q_{j-2})+\notag\\+
\sum_{j=0}^{n-2}a_{j,n-2}(q_{j-1}+q_{j-2})+\sum_{j=0}^{n-4}b_{j,n-2}(q_{j-1}+q_{j-2})+q_{n-1}+q_{n-2}.\label{pp}
\end{gather}
If the statement is proved for $k\le n-1$ then by \eqref{pp}, \eqref{p} we have
\begin{gather*}
\sum_{j=0}^na_{j,n}(q_{j-1}+q_{j-2})+\sum_{j=0}^{n-2}b_{j,n}(q_{j-1}+q_{j-2})=
4p_{n-1}+2p_{n-2}+q_{n-1}+q_{n-2}=2p_n.
\end{gather*}
This completes the proof.
\end{proof}
\end{Le}
\begin{Le}\label{sumS}
Let $0<y<q$ be an integer and $\lambda_i\in\{-1,1\}.$ Let $$0<x_1<x_2<\ldots<x_m$$ be a sequence of integers such that $x_i+y<x_{i+1}$ for $1\le i<m$, $x_m+y-x_1\le q$ then
\begin{gather*}
\left|\sum_{i=1}^m\lambda_iS(x_i,y)\right|\le y\sqrt{Bqm}.
\end{gather*}
\begin{proof}
Let $J_i=[x_i+1,x_i+y]$ for $1\le i\le m$ then $\bigcup J_i\subseteq[x_1+1,x_1+q]$ and $J_i\bigcap J_j=\emptyset$ for any $1\le i,\,j\le m$ Let
\begin{gather*}
\lambda_i(x)=
\left\{
              \begin{array}{ll}
                \lambda_i, & \hbox{if $x\in J_i$;} \\
                0, & \hbox{else.}
              \end{array}
\right.
\end{gather*}
By \eqref{S}, the Cauchy-Schwarz inequality and the assumption $f\in\textit{F}_{A,B}(q)$  we have
\begin{gather*}
\left|\sum_{i=1}^m\lambda_iS(x_i,y)\right|=\left|\sum_{n\in\bigcup J_i}\sum_{k=0}^{y-1}\lambda_i(n)f(n+k)\right|\le
\sqrt{\sum_{n\in\bigcup J_i}\lambda_i^2(n)}\sqrt{\sum_{n\in\bigcup J_i}\left|\sum_{k=0}^{y-1}f(n+k)\right|^2}\le\\\le
\sqrt{my}\sqrt{\sum_{n=1}^q\left|\sum_{k=0}^{y-1}f(n+k)\right|^2}\le y\sqrt{Bqm}.
\end{gather*}
\end{proof}
\end{Le}
To prove Theorem \ref{Th1} we must estimate $\frac{1}{K_0}\left|-T_{-}(a,K_0)+T_{+}(a+N,K_0)\right|$ (see \eqref{1estimate}). Put
\begin{gather*}
\Sigma_{-}^{(n)}(a,K_i)=\sum_{j=0}^na_{jn}S_{-}^{(j)}(a,K_i)-\sum_{j=0}^{n-2}b_{jn}S_{+}^{(j)}(a,K_i)\\
\Sigma_{+}^{(n)}(a,K_i)=\sum_{j=0}^na_{jn}S_{+}^{(j)}(a,K_i)-\sum_{j=0}^{n-2}b_{jn}S_{-}^{(j)}(a,K_i).
\end{gather*}
Applying Lemma \ref{recurs}  and formula \eqref{Tn}, we have
\begin{gather*}
T_{-}(a,K_0)=T_{-}^{(n)}(a,K_{1})=\alpha_nT_{-}^{(n)}(a,K_{2})-\beta_nT_{+}^{(n)}(a,K_{2})+\Sigma_{-}^{(n)}(a,K_1)=T_1+S_1
\end{gather*}
with
\begin{gather}\label{A1B1}
T_1=\alpha_nT_{-}^{(n)}(a,K_{2})-\beta_nT_{+}^{(n)}(a,K_{2}),\quad S_1=\Sigma_{-}^{(n)}(a,K_1).
\end{gather}
Consecutive application of Lemma \ref{recurs} will give us an upper bound for $T_{-}(a,K_0).$
If
\begin{gather*}
T_i=A_iT_{-}^{(n)}(a,K_{i+1})-B_iT_{+}^{(n)}(a,K_{i+1}), \quad 1\le i\le \tau-1,
\end{gather*}
then for $1\le i\le \tau-2$ by Lemma \ref{recurs} one has
\begin{gather*}
T_i=A_i\left(\alpha_nT_{-}^{(n)}(a,K_{i+2})-\beta_nT_{+}^{(n)}(a,K_{i+2})+\Sigma_{-}^{(n)}(a,K_{i+1})\right)-\\-
B_i\left(\alpha_nT_{+}^{(n)}(a,K_{i+2})-\beta_nT_{-}^{(n)}(a,K_{i+2})+\Sigma_{+}^{(n)}(a,K_{i+1})\right)=
\left(\alpha_nA_i+\beta_nB_i\right)T_{-}^{(n)}(a,K_{i+2})-\\-\left(\alpha_nB_i+\beta_nA_i\right)T_{+}^{(n)}(a,K_{i+2})+
A_i\Sigma_{-}^{(n)}(a,K_{i+1})-B_i\Sigma_{+}^{(n)}(a,K_{i+1})=T_{i+1}+S_{i+1}.
\end{gather*}
with
\begin{gather*}
S_{i+1}=A_i\Sigma_{-}^{(n)}(a,K_{i+1})-B_i\Sigma_{+}^{(n)}(a,K_{i+1}).
\end{gather*}
So
\begin{gather*}
\left\{
  \begin{array}{ll}
   A_{i+1}=\alpha_nA_i+\beta_nB_i , &  \\
   B_{i+1}=\alpha_nB_i+\beta_nA_i, &  \\
   A_1=\alpha_n,\;B_1=\beta_n.
    \end{array}
\right.
\end{gather*}
By \eqref{alphabeta} we have
\begin{gather}\label{AB}
A_i=\frac{q_n^i+1}{2},\quad B_i=\frac{q_n^i-1}{2}.
\end{gather}
Let
\begin{gather}\label{r}
r=\left\lceil\frac{\tau}{2}\right\rceil
\end{gather}
be the number of steps. As $\tau\ge3$ then $r\le \tau-1.$ So
\begin{gather*}
T_{-}(a,K_0)=T_r+\Sigma_{-}^{(n)}(a,K_1)+\sum_{i=2}^rS_i=\\=
A_rT_{-}^{(n)}(a,K_{r+1})-B_rT_{+}^{(n)}(a,K_{r+1})+\Sigma_{-}^{(n)}(a,K_1)+
\sum_{i=2}^r\left(A_{i-1}\Sigma_{-}^{(n)}(a,K_{i})-B_{i-1}\Sigma_{+}^{(n)}(a,K_{i})\right)
\end{gather*}
and
\begin{gather*}
T_{-}(a+N,K_0)=
A_rT_{-}^{(n)}(a+N,K_{r+1})-B_rT_{+}^{(n)}(a+N,K_{r+1})+\Sigma_{-}^{(n)}(a+N,K_1)+\\+
\sum_{i=2}^r\left(A_{i-1}\Sigma_{-}^{(n)}(a+N,K_{i})-B_{i-1}\Sigma_{+}^{(n)}(a+N,K_{i})\right).
\end{gather*}
So
\begin{gather}\label{t-t}
\left|T_{-}(a,K_0)-T_{-}(a+N,K_0)\right|\le\Sigma_1+\Sigma_2+\sum_{i=2}^r\Sigma_3(i)
\end{gather}
with
\begin{gather*}
\Sigma_1=\left|A_rT_{-}^{(n)}(a,K_{r+1})-B_rT_{+}^{(n)}(a,K_{r+1})-
A_rT_{-}^{(n)}(a+N,K_{r+1})+B_rT_{+}^{(n)}(a+N,K_{r+1})\right|,\\
\Sigma_2=\left|\Sigma_{-}^{(n)}(a,K_1)-\Sigma_{-}^{(n)}(a+N,K_1)\right|,\\
\Sigma_3(i)=\left|A_{i-1}\Sigma_{-}^{(n)}(a,K_{i})-B_{i-1}\Sigma_{+}^{(n)}(a,K_{i})-
A_{i-1}\Sigma_{-}^{(n)}(a+N,K_{i})+B_{i-1}\Sigma_{+}^{(n)}(a+N,K_{i})\right|.
\end{gather*}
Trivially we obtain by \eqref{Tn} and the definition of $T_{+}, T_{-}$ (see \eqref{T+}, \eqref{T-}) that
\begin{gather}\label{sum1}
\Sigma_1\le 2A\frac{K_r\left(K_r+1\right)}{2}\left(A_r+B_r\right)=AK_r\left(K_r+1\right)\left(A_r+B_r\right).
\end{gather}
By \eqref{Sn}, \eqref{cq}, Lemma \ref{sumS}, Lemma \ref{Pn}  we have
\begin{gather}
\Sigma_2\le
\sum_{j=0}^na_{j,n}\left|S_{-}^{(j)}(a,K_1)-S_{-}^{(j)}(a+N,K_1)\right|+
\sum_{j=0}^{n-2}b_{j,n}\left|S_{+}^{(j)}(a,K_1)-S_{+}^{(j)}(a+N,K_1)\right|\le\notag\\\le
\sqrt{2Bq}\left(K_1+\frac{1}{2}\right)
\Biggl(\sum_{j=0}^na_{j,n}(q_{j-1}+q_{j-2})+\sum_{j=0}^{n-2}b_{j,n}(q_{j-1}+q_{j-2})\Biggl)=\notag\\=
\sqrt{2Bq}\left(2K_1+1\right)p_n.\label{sum2}
\end{gather}
Applying \eqref{AB} we have $A_i=B_i+1$ so
\begin{gather*}
\Sigma_3(i)\le B_{i-1}\left|\Sigma_{-}^{(n)}(a,K_{i})-\Sigma_{+}^{(n)}(a,K_{i})-\Sigma_{-}^{(n)}(a+N,K_{i})+\Sigma_{+}^{(n)}(a+N,K_{i})\right|+\\+
\left|\Sigma_{-}^{(n)}(a,K_{i})-\Sigma_{-}^{(n)}(a+N,K_{i})\right|.
\end{gather*}
By \eqref{Sn}, \eqref{cq}, Lemma \ref{sumS}, Lemma \ref{Pn} we have
\begin{gather}\label{sum3}
\Sigma_3(i)\le
\sqrt{4Bq}B_{i-1}\left(2K_i+1\right)p_n+\sqrt{2Bq}\left(2K_i+1\right)p_n.
\end{gather}
Using \eqref{t-t}, \eqref{sum1}, \eqref{sum2}, \eqref{sum3} we have
\begin{gather}
\left|T_{-}(a,K_0)-T_{-}(a+N,K_0)\right|\le AK_r\left(K_r+1\right)\left(A_r+B_r\right)+\notag\\+
\sqrt{2Bq}p_n\sum_{i=1}^r\left(2K_i+1\right)+2\sqrt{Bq}p_n\sum_{i=2}^rB_{i-1}\left(2K_i+1\right)\label{T--T-}.
\end{gather}
By \eqref{Ki}, \eqref{AB} we have
\begin{gather*}
K_r\left(K_r+1\right)\left(A_r+B_r\right)=q_n^r\frac{q_n^{\tau-r}-1}{2}\frac{q_n^{\tau-r}+1}{2}=
\left(K_0-\frac{q_n^{r}-1}{2}\right)\frac{q_n^{\tau-r}+1}{2}=\\=K_0\frac{q_n^{\tau-r}+1}{2}-\frac{1}{2}K_0-\frac{q_n^r-q_n^{\tau-r}}{4}=
K_0\frac{q_n^{\tau-r}}{2}-\frac{q_n^r-q_n^{\tau-r}}{4}.
\end{gather*}
And
\begin{gather*}
B_{i-1}\left(2K_i+1\right)=\frac{q_n^{i-1}-1}{2}q_n^{\tau-i}=\frac{q_n^{\tau}-1-q_n^{\tau-i+1}+1}{2q_n}=\frac{2K_0+1}{2q_n}-
\frac{q_n^{\tau}}{2q_n^i}.
\end{gather*}
So we obtain
\begin{gather*}
\left|T_{-}(a,K_0)-T_{-}(a+N,K_0)\right|\le
A\left(K_0\frac{q_n^{\tau-r}}{2}-\frac{q_n^r-q_n^{\tau-r}}{4}\right)+2\sqrt{Bq}p_n\frac{2K_0+1}{2q_n}(r-1)+\\+
\sqrt{Bq}p_nq_n^{\tau}\left(\sqrt{2}\sum_{i=1}^r\frac{1}{q_n^i}-\sum_{i=2}^r\frac{1}{q_n^i}\right).
\end{gather*}
Applying \eqref{Ki} we have
\begin{gather*}
\sqrt{Bq}p_nq_n^{\tau}\left(\sqrt{2}\sum_{i=1}^r\frac{1}{q_n^i}-\sum_{i=2}^r\frac{1}{q_n^i}\right)\le
\sqrt{Bq}\left(2K_0+1\right)\left(\sqrt{2}-1\right)\frac{p_n}{q_n-1}+\sqrt{Bq}\left(2K_0+1\right)\frac{p_n}{q_n}.
\end{gather*}
So by \eqref{r} we get
\begin{gather*}
\left|T_{-}(a,K_0)-T_{-}(a+N,K_0)\right|\le
AK_0\frac{q_n^{\tau-r}}{2}+
\sqrt{Bq}\left(2K_0+1\right)\left(\frac{p_n}{q_n}r+\left(\sqrt{2}-1\right)\frac{p_n}{q_n-1}\right).
\end{gather*}
Using \eqref{1estimate} we have
\begin{gather}\label{2estimate}
\left|\sum_{n=a+1}^{a+N}f(n)\right|\le\sqrt{B q}\sqrt{\frac{N}{K_0}}+
A\frac{q_n^{\tau-r}}{2}+
\sqrt{Bq}\left(2+\frac{1}{K_0}\right)\left(\frac{p_n}{q_n}r+\left(\sqrt{2}-1\right)\frac{p_n}{q_n-1}\right).
\end{gather}
By \eqref{r}, \eqref{NN}, \eqref{Nq2} we get
$$
q_n^{\tau-r}=q_n^{\tau-\left\lceil\frac{\tau}{2}\right\rceil}\le\sqrt{q_n^{\tau}}\le\sqrt{q}.
$$
So by \eqref{NK0} we obtain
\begin{gather}\label{3estimate}
\left|\sum_{n=a+1}^{a+N}f(n)\right|\le\sqrt{B q}\sqrt{q_n+\frac{1}{q_n^{2}}}+
A\frac{\sqrt{q}}{2}+
\sqrt{Bq}\left(2+\frac{q_n+\frac{1}{q_n^{2}}}{N}\right)\left(\frac{p_n}{q_n}r+\left(\sqrt{2}-1\right)\frac{p_n}{q_n-1}\right).
\end{gather}
By the definition of $\tau, r$ we have
$$r\le\frac{\log q}{2\log q_n}.$$
So
\begin{gather}
\left|\sum_{n=a+1}^{a+N}f(n)\right|\le
\sqrt{B}\delta_n\sqrt{q}\log q+\sqrt{B q}\sqrt{q_n+\frac{1}{q_n^{2}}}+A\frac{\sqrt{q}}{2}+
2\sqrt{Bq}\left(\sqrt{2}-1\right)\frac{p_n}{q_n-1}+\notag\\+
\sqrt{Bq}\frac{1}{N}\left(q_n+\frac{1}{q_n^{2}}\right)\left(\frac{p_n}{q_n}\frac{\log q}{2\log q_n}+\left(\sqrt{2}-1\right)\frac{p_n}{q_n-1}\right).
\label{4estimate}
\end{gather}
Applying \eqref{N}, we have
\begin{gather*}
\left|\sum_{m=a+1}^{a+N}f(m)\right|\le\sqrt{B}\delta_n\sqrt{q}\log q+
\left(\sqrt{B}\sqrt{q_n+\frac{1}{q_n^{2}}}+\sqrt{B}(\sqrt{2}-1)\frac{2p_n}{q_n-1}+\frac{1}{2}A\right)\sqrt{q}+\\+,
2\sqrt{B}\left(q_n+\frac{1}{q_n^{2}}\right)\left(\frac{p_n}{q_n}\frac{\log q}{2\log q_n}+\left(\sqrt{2}-1\right)\frac{p_n}{q_n-1}\right).
\end{gather*}
Theorem \ref{Th1} is proved.
\section{On the constant in the P\'{o}lya--Vinogradov inequality}
Burgess \cite{Burgess} proved that for a nonprincipal character $\chi \pmod{q}$ one has $\chi\in\textit{F}_{1,1}(q).$  Applying this result to \eqref{DW} and \eqref{BR}, we have
\begin{gather}\label{DW1}
\left|\sum_{n=a+1}^{a+N}\chi(n)\right|\le\frac{1}{2\log2}\sqrt{q}\log q+3\sqrt{q},
\end{gather}
\begin{gather}\label{BR1}
\left|\sum_{n=a+1}^{a+N}\chi(n)\right|\le\frac{1}{3\log3}\sqrt{q}\log q+6.5\sqrt{q}.
\end{gather}
But this result is not the best one. Let
$$S_{\chi}=\max_{0\le M<N\le q}\left|\sum_{n=M}^N\chi(n)\right|, \qquad T_{\chi}=\max_{N}\left|\sum_{a=0}^N\chi(a)\right|.$$
Granville and Soundararajan \cite{Granville-Sound} obtained two inequalities
\begin{gather}\label{GS1}
T_{\chi}\le\left(\frac{69}{70}\frac{c}{\pi\sqrt{3}}+o(1)\right)\sqrt{q}\log q \quad \mbox{if}\quad \chi(-1)=1,
\end{gather}
and
\begin{gather}\label{GS2}
T_{\chi}\le\left(\frac{c}{\pi}+o(1)\right)\sqrt{q}\log q \quad \mbox{if}\quad \chi(-1)=-1,
\end{gather}
where
$$
c=
\left\{
              \begin{array}{ll}
                \frac{1}{4}, & \hbox{if q is a cubefree;} \\
                \frac{1}{3}, & \hbox{else.}
               \end{array}
\right.
$$
Up to now this result is the best-known one. Pomerance proved (see~\cite{Pomerance}) numerically explicit version of the P\'{o}lya--Vinogradov inequality
$$
S_{\chi}\le\frac{2}{\pi^2}\sqrt{q}\log q+\frac{4}{\pi^2}\sqrt{q}\log\log q+\frac{3}{2}\sqrt{q}\quad \mbox{if}\quad \chi(-1)=1
$$
and
$$
S_{\chi}\le\frac{1}{2\pi}\sqrt{q}\log q+\frac{1}{\pi}\sqrt{q}\log\log q+\sqrt{q}\quad \mbox{if}\quad \chi(-1)=-1.
$$
Up to now these bounds are the best-known numerically explicit versions of the P\'{o}lya--Vinogradov inequality.
These inequalities are weaker then \eqref{GS1},\,\eqref{GS2} but better then \eqref{BR1}.
Applying Theorem \ref{Th1} we improve \eqref{BR1}.
\begin{Sl}\label{Sl1}
For any $n\ge0$, $q\ge q_n^6$ and any nonprincipal character $\chi \pmod{q}$, we have
\begin{gather*}
\left|\sum_{m=a+1}^{a+N}\chi(m)\right|\le\delta_n\sqrt{q}\log q+
\left(\sqrt{q_n+\frac{1}{q_n^{2}}}+(\sqrt{2}-1)\frac{2p_n}{q_n-1}+\frac{1}{2}\right)\sqrt{q}+\psi_n(q,1),
\end{gather*}
\end{Sl}


\textit{D.A. Frolenkov}\\
\textit{Department of Number theory}\\
\textit{Moscow State University}\\
\textit{e-mail: frolenkov\underline{ }adv@mail.ru}

\end{document}